\documentclass{elsart}
\usepackage{amssymb}

\def\reff#1{(\ref{#1})}
\def\basis{{\rm Basis}\,}
\def\life{{\rm Life}\,}
\def\birth{{\rm Birth}\,}
\def\death{{\rm Death}\,}
\def\X{{\mathcal X}}
\def\XX{{\bf X}}
\def\vep{{\varepsilon}}
\def\E{{\mathbb E}}
\def\P{{\mathbb P}}
\def\R{{\mathbb R}}
\def\Z{{\mathbb Z}}

\def\N{{\mathbb N}}

\def\Y{{\bf Y}}
\def\C{{\bf C}}
\def\D{{\bf D}}
\def\B{{\bf B}}
\def\G{{\bf G}}
\def\K{{\bf K}}
\def\L{{\mathcal L}}
\def\one{{\bf 1}\hskip-.5mm}
\def\be{{\beta}}
\def\de{{\delta}}
\def\la{{\lambda}}
\def\ga{{\gamma}}
\def\ka{{\kappa}}
\def\th{{\theta}}
\def\rate{{\e^{- \beta|\ga|}}}
\def\A{{\bf A}}
\def\B{{\bf B}}
\def\C{{\bf C}}
\def\D{{\bf D}}

\begin{document}

\runauthor{Ferrari and Picco}
\begin{frontmatter}
\title{Poisson approximation for large-contours 
\\ in low-temperature Ising models}

\author[SP]{Pablo A. Ferrari} 
\author[MRS]{Pierre Picco}
\address[SP]{Universidade de S\~{a}o Paulo, IME, {Cx.} Postal 66281,
  05315-970, S\~ao Paulo, Brazil; email: 
  pablo@ime.usp.br, http://www.ime.usp.br/\~{}pablo} 
\address[MRS]{CNRS-CPT Luminy Case 907, 13288 Marseille, France}

\begin{abstract} We consider the contour representation of the
  infinite volume Ising model at any fixed inverse temperature
  $\beta>\beta^*$, the solution of $ \sum_{\theta : \theta\ni 0}
   \e^{-\beta|\theta|}=1$. Let $\mu$ be the infinite-volume ``$+$'' measure.
  Fix $V\subset\Z^d$, $\lambda>0$ and a (large) $N$ such that calling
  $\G_{N,V}$ the set of contours of length at least $N$ intersecting
  $V$, there are in average $\lambda$ contours in $\G_{N,V}$ under
  $\mu$. We show that the total variation distance between the law of
  $(\ga:\ga\in\G_{N,V})$ under $\mu$ and a Poisson process is bounded
  by a constant depending on $\beta$ and $\lambda$ times $\e^{-(\beta
    -\beta^*)N}$.  The proof builds on the Chen-Stein method as
  presented by Arratia, Goldstein and Gordon. The control of the
  correlations is obtained through the loss-network space-time
  representation of contours due to Fern\'andez, Ferrari and Garcia.
\end{abstract}


\end{frontmatter}

\section{Introduction}

The infinite volume Ising model is one of the most studied in
Statistical Physics.  One of the fundamental historical fact was the
Peierls argument for the existence of more than one phase at low
enough temperature ---which we assume throughout this paper--- in any
dimension $d \ge 2$. This argument was set by Dobrushin as the
existence of at least two extremal Gibbs states, $\nu^+$ and $\nu^-$
that are obtained as infinite volume limits of Gibbs measures with $+$
resp. $-$ boundary conditions. In the Physics literature, a
configuration of spins which is typical with respect to $\nu^+$, is
frequently described as a set of islands of $-$ within a sea of $+$ ;
the contours being the boundaries of the islands.  Various questions
can be asked about typical configuration of contours at low
temperature.  A very simple one that does not seem to have attracted
attention is: what is the number of occurrences of contours larger
than a fixed (big) $N$, that intersect a fixed volume $V=V(N)$ say
centered at the origin. To describe the physical idea of rarity of the
island of $-$ , a Poisson type behavior naturally comes in mind.  We
prove indeed that the distribution of these occurrences,
under the condition that the mean number of occurrences in the volume
$V(N)$ is a fixed value $\lambda$, approaches sharply a Poisson process
of mean $\lambda$.  The paper is organized as follows: in section
2, we give the definitions, state the main result and the
Arratia-Goldstein-Gordon formulation of the Chen-Stein Method used in
the proof.  In section 3 we recall the graphical construction of the
loss-network space-time representation of contours of
Fern\'andez-Ferrari-Garcia which we apply to obtain the relevant
bounds. In section 4 we prove the main result.

\section{ Definitions and the main result}

Peierls introduced a map between typical configurations of $\nu^+$
or $\nu^-$ into an ensemble of objects ---the contours--- interacting
only by perimeter-exclusion.  See, for instance, Section 5B of
Dobrushin \cite{dob96}, for a concise and rigorous account of this
mapping.  Contours are hyper surfaces formed by a finite number of
$(d-1)$-dimensional unit cubes ---\emph{links} for $d=2$,
\emph{plaquettes} for higher dimensions--- centered at points of
$\Z^{d}$ and perpendicular to the edges of the dual lattice
$\Z^d+({1\over 2},\cdots,{1\over 2})$.  To formalize their definition,
let us call two plaquettes \emph{adjacent} if they share a
$(d-2)$-dimensional face.  A set of plaquettes, $\gamma$, is
\emph{connected} if for any two plaquettes in $\gamma$ there exists a
sequence of adjacent plaquettes in $\gamma$ joining them.  The set
$\gamma$ is \emph{closed} if every $(d-2)$-dimensional face is covered
by an even number of plaquettes in $\gamma$.  \emph{Contours} are
connected and closed sets of plaquettes.  For example, in two
dimensions contours are closed polygonals.  Two contours $\ga$ and
$\theta$ are said \emph{compatible} if no plaquette of $\ga$ is
adjacent to a plaquette of $\theta$. In this case we write
$\ga\sim\th$. In two dimensions, therefore, contours are compatible if
and only if they do not share the endpoint of a link.  In three
dimensions two compatible contours can share vertices, but not sides
of plaquettes.  Ising spin configurations in a bounded region with
``+'' (or ``$-$'') boundary condition are in one-to-one correspondence
with families of pairwise compatible contours.

The set of possible
contours contained in $\Lambda\subset\Z^d$ will be denoted
$\G[\Lambda]$. The set of configurations of compatible
contours is  
$ \X[\Lambda] = \{ \eta \in \{0,1\}^{\G[\Lambda]}\,;\,
  \eta(\ga)\,\eta(\theta) = 
  0 \mbox{ if } \ga \not\sim\theta\}$. 
We denote $|\ga|$ the number of plaquettes of the contour $\ga$ and
define the finite-volume Gibbs measure $\mu_\Lambda$ on
$\X[\Lambda]$: for $\eta\in\X[\Lambda]$,
\begin{equation}
  \label{141}
  \mu_\Lambda(\eta) = {\exp\Bigl(-\beta \sum_{\ga} 
    |\ga|\,\eta(\ga)\Bigr)\over Z_\Lambda} 
\end{equation}
We fix $\beta>\beta^*$, the solution of $\alpha_0(\be)=1$, where
$\alpha_0(\beta):= \sum_{\theta : \theta\ni 0} 
\,\e^{-\beta|\theta|}\,$.
Then there exists a unique weak limit $\mu=\lim_{\Lambda\to\Z^d}
\mu_\Lambda$ \cite{ffg,ffg2}. Notice that $\be^*$ is strictly bigger
than $\be_P$, the Peierls inverse-temperature, which is the infimum of
the $\be$ making the previous sum 
 finite. Let $p_\ga = \int
\mu(d\eta)\,\eta(\ga)$ be the probability of the presence of contour
$\ga$ under $\mu$. Fix $N>0$, let $V=V(N)\subset\Z^d$ and $\lambda>0$
be such that if we define
$\G\,=\,\G(N,V,\lambda)\, 
:=\, \{\ka \in \G[\Z^d]\,:\, |\ka|\ge N\hbox{ and }\ka\cap V\neq\emptyset\}$
as the set of contours with length at least $N$ intersecting $V$, then
$ \sum_{\ga\in\G} p_\ga=\lambda$. 
Hence $V$ is a set such that in average there are $\lambda$ contours
of length at least $N$ intersecting $V$. 
Define the processes 
 $ \XX:=(\eta(\ga)\,:\,\ga\in \G),\;\Y:=(Y(\ga)\,:\,\ga\in \G)$
where $\eta$ has distribution $\mu$ and $Y(\ga)$ are iid with Poisson
distribution of mean $p_\ga$.
We denote $\L({\bf X})$ the law of a process ${\bf X}$ and $\|
\L(\XX)-\L(\Y)\|_{\rm TV}$ the total variation distance between the
laws of the processes $\XX$ and $\Y$ \cite{agg}. 

\begin{thm}
      \label{191} Fix $\beta>\beta^*$, then 
for any $\be'\in(\be^*,\be)$, if  $M(\be,\be')$ is given in
\reff{mbe} below we have:

\begin{equation}
  \label{125}
  \| \L(\XX)-\L(\Y)\|_{\rm TV} \; 
\leq \; M(\beta,\be'){N^{d+1}}{\lambda}\,
  \e^{-(\be-\be')N} 
\end{equation}
\end{thm}

The proof of this theorem is an application of the Chen--Stein method
\cite{stein,chen,agg,bhj} as proposed in \cite{agg}. In order to
describe it, let the distance between contours be
$d(\ga,\th):=\min\{|x-y|\,:\,x\in\ga,\,y\in\th\}$. For $D=D(N)>0$ to
be fixed later we define $\B_\ga\,=\,\B_\ga(N,D) = \{\th\in\G\,:\,
d(\ga,\th) < D\}$, a ``$D$-neighborhood'' of $\ga$.  Let $p_{\ga\th} =
\int\mu(d\eta)\eta(\ga)\eta(\th)$ and define $b_i=b_i(N,D)$ by
\begin{eqnarray}
\nonumber
b_1 &:=& \sum_{\ga\in \G}\; \sum_{\th\in \B_\ga} p_\ga p_\th
\qquad\qquad\qquad\qquad 
b_2 \;:=\; \sum_{\ga\in \G}\, \sum_{\ga\neq\th\in \B_\ga} p_{\ga\th} \\
b_3 &:=& \sum_{\ga\in \G} \E\Bigl| \E\{\eta(\ga) - p_\ga\;|\; 
\sigma(\eta(\th)\,:\, \th \in \G\setminus \B_\ga)\}\Bigr|
\nonumber
\end{eqnarray}

Theorem 2 of \cite{agg} shows that the total variation distance
between the Poisson process $\Y$ and the process $\XX$ is dominated by
$2(2b_1+2b_2+b_3)$. To prove Theorem \ref{191} we choose $D$ equal to
$N$ times an appropriate constant and then find bounds on $b_i$. The
bound on $b_1$ does not offer problems, but the bound on $b_2$
requires the inequality $p_{\ga\th} \le p_\ga \e^{-\be|\th|}+p_\th
\e^{-\be|\ga|}$ proved in Lemma \ref{pab} below. The Chen-Stein method
has been applied mostly to examples where $b_3=0$ \cite{agg}, while in
our case $b_3\neq 0$.  To dominate $b_3$ we control the correlations
between ``far'' contours in Lemma \ref{b3} below. The proof of those
lemmas is based on the graphical representation of a loss network, a
Markov process having $\mu$ as invariant measure \cite{ffg,ffg2}.
There are many papers dealing with Poisson approximations; we limit
ourselves to quote the books \cite{a} and \cite{bhj} ---with a nice
introduction explaining the Chen-Stein method---, the original paper
of Chen \cite{chen} and a monograph of Stein \cite{stein}.
There are a few results for random fields. The number of ``$-$''
spins in the ``$+$'' measure of the Ising model at low temperature
and/or high external ``$+$'' magnetic field converges to Poisson
\cite{bg,ghosu}. In this case the variables are positively associated
and the Chen-Stein method works via cluster expansions giving explicit
convergence rates \cite{ghosu}. In \cite{ffg2} the approximation is
established for the contours of the Ising model at low temperature and
without magnetic field. The available bounds for the total variation
distance are of the order of the inverse of the volume of the observed
set $V$ \cite{ghosu,ffg2}.

\section{Graphical construction}

In this section we quickly review the construction and some results of
\cite{ffg} relevant to our proof. To each contour $\ga \in \G(\Z^d)$
we associate an independent (of everything) marked Poisson process
$N_{\ga}$ on $\R$ with rate $\e^{-\be |\ga|}$. We call
$T_{k}(\ga)\in\R$, $\ga\in\G(\Z^d)$, the ordered time-events of
$N_\ga$ with the convention that $T_{0}(\ga)<0<T_{1}(\ga)$.  For each
occurrence time $T_{i}(\ga)$ of the process $N_{\ga}$ we choose an
independent mark $S_{i}(\ga)$ exponentially distributed with mean 1.
At the Poisson time-event $T_i(\ga)$ a contour $\ga$ appears and it
lasts $S_i(\ga)$ time units.  The resultant object is the random
family $\C=\bigl\{\{(\ga, T_{i}(\ga),S_{i}(\ga)): i \in \Z \}:\ga\in
\G\bigr\}$. A marked point $(\ga,T_k(\ga),S_k(\ga))\in\C$ is
identified with $\ga \times [T_k(\ga), T_k(\ga)+S_k(\ga)]$, the
\emph{cylinder} with \emph{basis} $\ga$, \emph{birth-time} $T_k(\ga)$
and \emph{lifetime} $S_k(\ga)$.  The \emph{life} of the cylinder is
the time interval $[T_k(\ga),T_k(\ga)+S_k(\ga)]$. For a generic
cylinder $C=(\ga,t,s)$, we use the notation $\basis(C) = \ga$,
$\birth(C)= t$, $\death(C)= t+s$, $\life(C) = [t,t+s]$.
In the sequel $\P$ and $\E$ are the probability and expectation in the
space where $\C$ is defined.

For $t\in\R$ we define 
 $ \xi_t(\ga,\C) = \sum_{C\,\in\,\C} \one\{\basis(C)=\ga,
\life(C)\ni t\}.$ 
The above process, called the \emph{free network}, is a product of
independent stationary birth-and-death processes on $\N^{\G[\Z^d]}$ with
$t\in\R$ whose generator is given by
 \begin{eqnarray}
 \nonumber
A^{0} f(\xi) &=& \sum_{\ga \in \G[\Z^d]} \rate [f(\xi+\de_\ga) -
f(\xi)]
\; +\; \sum_{\ga \in \G[\Z^d]} \xi(\ga)[f(\xi-\de_\ga) -
f(\xi)].\label{14b} 
\end{eqnarray}
The invariant (and reversible) measure for this process is the product measure
$\mu^{0}$ on $\N^{\G[\Z^d]}$ with Poisson marginals:
$\mu^{0}\bigl\{\xi(\ga)=k\bigr\} \,=\,\exp\left(\rate\right)(\rate)^{k}/k!$ 
In particular, for any $t\in\R$, $\E f(\xi_t) = \mu^0f$.

We say that cylinders $C$ and $C'$ are \emph{incompatible} and write
$C'\not\sim C$ if and only if $\basis(C)\not\sim\basis(C')$ and
$\life(C)\cap\life(C')\neq\emptyset$. We say that two sets of
cylinders $\A$ and $\A'$ are \emph{incompatible} and write
$\A\not\sim\A'$ if there is a cylinder in $\A$ incompatible with a
cylinder in $\A'$. Otherwise we use the sign $\sim$ for compatibility.
For any cylinder $C$ define the set of {\it
  ancestors} of $C$ as the set of cylinders in $\C$ born before $C$
that are incompatible with $C$: ${\A}^{C}_1 = \{C'\in\C\, ;\,
C'\not\sim C, \birth(C')<\birth(C)\}$.  Recursively for $n\ge 2$, the
$n$th generation of ancestors of $C$ is ${\A}_n^C = \{ C'':C''\in{\bf
  A}^{C'}_1\hbox{ for some } C' \in {\bf A}_{n-1}^C \}$.  Let the
\emph{clan} of $C$ be the union of its ancestors: $\A^C=\bigcup_{n\ge
  1}{\A}_{n}^C$. 
Under the condition $\be>\be^*$ all cylinders in $\C$ have a finite
clan with probability one. This property called \emph{no backwards
  oriented percolation} is essential to show that the loss network
$\eta_t$ can be constructed in a stationary way for $t\in\R$. 
Assume that there is no backwards oriented percolation.
The construction is as follows. All cylinders in $\C$ are classified
as \emph{kept} and \emph{erased}. Since all clans are finite, we can
write $\C=\cup_{n\ge 0} \C_n$, where $\C_n:= \{C\in\C: \A_{n}^C \neq
\emptyset,\;\A_{n+1}^C =\emptyset\} $. Inductively we set $\K_0=\C_0$,
$\D_0=\emptyset$, $\K_n =\{C\in \C_n\setminus
\cup_{i=0}^{n-1}(\D_i\cup\K_i)\,:\, \{C\}\sim \cup_{i=0}^{n-1}\K_i\}$
and $\D_n =\C_n\setminus[\K_n\cup\, \cup_{i=0}^{n-1}(\D_i\cup\K_i)] $.
Let the set of kept cylinders be $\K=\cup_n \K_n$ and the set of
erased cylinders be $\D=\cup_n\D_n$. Clearly $\K\dot\cup\D=\C$. The
event $\{C\in\K\}$ is measurable with respect to the sigma field
generated by $\A^C$. In words, it is sufficient to know the (finite)
clan of $C$ to know if $C$ is kept or erased.
The stationary loss network is defined by
$\eta_t(\ga,\C) \; = \;\sum_{C\in\K} \one
  \Bigl\{\basis(C)=\ga,\, \life(C)\ni t\Bigr\}$.
The process $\eta_t$ is Markovian with generator
\begin{eqnarray}
  \label{144b}
A f(\eta) &=& \sum_{\ga \in \G} \rate\, \one\{\eta+\de_\ga\in\X\}\, 
[f(\eta+\de_\ga) - f(\eta)] \nonumber\\
&&\qquad\qquad\qquad\qquad +\; \sum_{\ga \in \G} \eta(\ga)
[f(\eta-\de_\ga) - f(\xi)].
\end{eqnarray}
The unique invariant (and reversible) measure for this process is the
Gibbs measure $\mu$.  In particular, for any $t\in\R$, $\E f(\eta_t) =
\mu f$. We study properties of $\mu$ by studying the law of $\eta_0$,
the stationary loss network at time zero.

The presence/absence of contours intersecting a region $\Lambda$ at
time $t$ depends only on the set $\A^{\Lambda,t}$, the union of the
clans of the cylinders of $\C$ with basis intersecting $\Lambda$ and
life containing time $t$, that is $\A^{\Lambda,t}:= \{C'\in\A^C\,:\,
\basis(C)\cap \Lambda \neq\emptyset,\, \life(C)\ni t\}$. In particular
$\eta_t(\ga)$ is a (deterministic) function of $\A^{\ga,t}$ defined by
 $\eta_t(\ga,\C)=\eta_t(\ga,\A^{\ga,t})$. The reason is
that in order to determine if $C\in\K$ it suffices to look at $\A^C$.
When $t=0$ we will use the notation $\A^\Lambda$ instead of
$\A^{\Lambda,0}$.

Let $C_{\ga,t,r} = (\ga,-t,t+r)$ be a cylinder with
$\birth(C_{\ga,t,r})=-t$, $\basis(C_{\ga,t,r})=\ga$,
$\life(C_{\ga,t,r})=t+r$. Conditioning to the birth time of the
(kept) cylinder alive at time zero with basis $\ga$, 
\begin{eqnarray}
p_\ga\;=\;
\int_{(\R^+)^2} \P(C_{\ga,t,r}\in \K\,|\, C_{\ga,t,r}\in \C)
\,\e^{-t}\,\e^{-r} \e^{-\be|\ga|}\,dt\,dr \label{rm2}
\end{eqnarray}

\noindent{\bf Coupling of clans. }
For disjoint sets $\Lambda$ and $\Upsilon$ of $\Z^d$ it is possible
to construct
$(\A^\Lambda,\A^\Upsilon,\widehat\A^\Lambda,\widehat\A^\Upsilon)$, a
coupling between four sets of cylinders satisfying
(a) $\A^\Lambda = \widehat\A^\Lambda$ and $\A^\Upsilon =
  \widehat\A^\Upsilon$, in distribution; 
(b) $\A^\Lambda\cup\A^\Upsilon = \A^{\Lambda\cup\Upsilon}$; 
(c) $\widehat\A^\Lambda$ and $\widehat\A^\Upsilon$ are
  \emph{independent}; 
(d) If $\widehat\A^\Lambda\sim\widehat\A^\Upsilon$, then the
  marginals \emph{coincide}: $
\widehat\A^\Lambda =\A^\Lambda \hbox{ and }\widehat\A^\Upsilon
=\A^\Upsilon  $.
We write $\P$
and $\E$ for the probability and the expectation of the coupling.
The following bound for the probability of incompatibility between
clans follows as in the proof of (2.13) of \cite{ffg2}: for any
$\be'\in(\be^*,\be)$, 
\begin{equation}
  \label{com}
  \P(\widehat\A^\Lambda \not\sim\widehat\A^{\Upsilon})\; \le
  \;2\,(1-\alpha_0(\be'))^{-2}\,\sum_{x\in\Lambda}\sum_{y\in\Upsilon}
 |x-y|\,
  \e^{-(\beta-\beta')|x-y|}\,.
\end{equation}

\section{ Proof of the theorem}

The proof of the theorem is based on a sequence of lemmata.
\begin{lem}
  \label{p10}
Denote $\rho=\alpha_0(\be)$. Then for $\ga$ such that
$|\ga|\ge N$, 
\begin{equation}
  \label{p11}
 \exp\{-(\be +\rho)
  |\ga|\}\;\le\; p_\ga\;\le\; 
  \exp\{-\be\, |\ga|\} 
\end{equation}
\end{lem}

\begin{pf}
  Bounding the conditional probability inside the integral in
  \reff{rm2} by one we get the upper bound. For the lower bound notice
  that the event $\{C_{\ga,t,r}\in \K\}$ is measurable with respect to
  the sigma field generated by the cylinders born before $-t$, the
  birth time of $C_{\ga,t,r}$. Hence
\begin{eqnarray}
    \label{rm7}
   \lefteqn{ \P(C_{\ga,t,r}\in \K\,|\, C_{\ga,t,r}\in \C) \;=\;
    \P\Bigl(
    \bigcap_{\th:\th\not\sim\ga}\{\eta_{-t}(\th,\A^{\th,-t})=0\}\Bigr) }\\
&&\ge \P\Bigl(
    \bigcap_{\th:\th\not\sim\ga}\{\xi_{-t}(\th,\C)=0\}\Bigr)
\;=\; \exp \Bigl(-\sum_{\th:\th\not\sim\ga} \e^{-\be|\th|}\Bigr) \;\ge\;
\e^{-\rho|\ga|}   
\label{ph4} 
\end{eqnarray}
The first inequality follows from $\eta_t(\ga)\le\xi_t(\ga)$; the
second identity from the distribution of $\xi_t$ described in Section 3.
Inserting this inequality in \reff{rm2} we get the left inequality in
\reff{p11}.
\end{pf}

\begin{lem}
\label{p38}
The following inequalities hold
  \begin{equation}
  \label{100}
|V| \,\le\,\la  \Bigl(\sum_{\ga\in\G_0} |\ga|^{-1}\,p_\ga\Bigr)^{-1}
\,\le\, \la  N\,\e^{(\be +\rho) N} 
\end{equation}
\end{lem}

\begin{pf}
  The first inequality is immediate. Bounding below the sum in
  \reff{100} by one term (choose any $\ga_0$ such that $|\ga_0|=N$) we
  dominate the middle term in \reff{100} by $|\ga_0|/p_{\ga_0}$ which
  using the left inequality in \reff{p11} is bounded by the rhs of
  \reff{100}.  
\end{pf}

\begin{lem}
\label{pab}
We have 
$p_{\ga\th}\,\le\, p_\ga\,\e^{-\be|\th|}+ p_\th\,\e^{-\be|\ga|}$
\end{lem}

\begin{pf}
 Notice that $p_{\ga\th} = 0$ if
$\ga\not\sim\th$.
Recall the notation used in \reff{rm2}. Consider compatible $\ga$ and
 $\th$ and condition to the birth times of the kept cylinders alive at
 time zero with bases $\ga$ and $\th$ to obtain as in \reff{rm2},
\begin{eqnarray}
p_{\ga\th}&=&\int_{(\R^+)^4} \P(C_{\ga,t,r}\in \K,\, C_{\th,s,w}\in\K\,|\,
C_{\ga,t,r}\in \C,\, 
C_{\th,s,w}\in\C)\, \nonumber\\
&&\qquad\qquad\qquad\qquad\times\,\e^{-t}\,\e^{-s}\,\e^{-r}\,\e^{-w}
\e^{-\be|\ga|}\,\e^{-\be|\th|}\,ds\,dt\,dr\,dw  \label{rm1}
\end{eqnarray}
For $s<t$ the event $\{C_{\ga,t,r}\in \K\}$ is measurable with respect
to the sigma field generated by $\{C\,:\,\birth(C)\le -t\}$. Hence it is
independent of $\{C_{\th,s,w}\in\C\}$ and we have 
\begin{eqnarray}
  \label{smt}
  \lefteqn{\one\{s<t\}\,\P(C_{\ga,t,r}\in \K,\, C_{\th,s,w}\in\K\,|\,
    C_{\ga,t,r}\in \C,\, C_{\th,s,w}\in\C)\nonumber}\\
&&\quad\le\; \one\{s<t\}\,\P(C_{\ga,t,r}\in \K\,|\,C_{\ga,t,r}\in\C,\,
    C_{\th,s,w}\in\C)\nonumber\\ 
&&\quad=\; \one\{s<t\}\,\P(C_{\ga,t,r}\in \K\,|\,C_{\ga,t,r}\in\C)\nonumber\\ 
&&\quad\le\;  \P(C_{\ga,t,r}\in \K\,|\,C_{\ga,t,r}\in\C)
\end{eqnarray}
Analogously $\one\{s\ge t\}\,\P(C_{\ga,t,r}\in \K,\,
C_{\th,s,w}\in\K\,|\, C_{\ga,t,r}\in \C,\,
C_{\th,s,w}\in\C)\,\le\,\P(C_{\th,s,w}\in \K\,|\,C_{\th,s,w}\in\C)$.
This and \reff{smt} allow us to factorize the integrals in \reff{rm1}
and then use \reff{rm2} to get the lemma.  
\end{pf}
\begin{lem} 
  \label{1b2} 
The following bound holds for $b_1$, for all $\vep\,\le\,d\,(\be+\rho)$
  \begin{eqnarray}
       b_1&\le&  \sum_{\ga\in \G}
(D+|\ga|)^d\, p_\ga\; \sum_{\th\in\G_0}
\,\e^{-\be|\th|}\label{bb1}\\
&\le& 2\, \lambda N\,\Bigl({\be+\rho\over
    \be-\be^*}\Bigr)\, \Bigl({D+1\over \vep}\Bigr)^d
    \,\e^{-(\be-\be^*-\vep)N } \label{2b2}
  \end{eqnarray}
\end{lem}

\begin{pf} Inequality \reff{bb1} follows immediately from \reff{p11}.
To prove \reff{2b2} it is sufficient to show that for all
$\vep\,<\,d\,(\be+\rho)$,
  \begin{equation}
    \label{ph6}
    \sum_{\ga\in \G} 
(D+|\ga|)^d\, p_\ga\;\le\; 2\, \lambda N\,\Bigl({\be+\rho\over
    \be-\be^*}\Bigr)\, \Bigl({D+1\over \vep}\Bigr)^d\, \exp(\vep N)\,. 
  \end{equation}

Using  $(D+|\ga|)^d\,\le\, (D+1)^d\; \e^{|\ga|\vep d}/\vep^d$ and \reff{100},
the lhs of \reff{ph6} is bounded by
\begin{equation}
  \label{ph9}
  {(D+1)^d\over \vep^d}\; \sum_{\ga\in \G}\e^{|\ga|\vep d}\, p_\ga \;\le\;
  \lambda {(D+1)^d\over {\vep^d}}\;
{\displaystyle{\sum_{\ga\in\G_0}\e^{|\ga|\vep 
      d}\, p_\ga} \over \displaystyle{\sum_{\ga\in\G_0}\,
    |\ga|^{-1}\,p_\ga} } 
\end{equation}

Divide the sum in the numerator in two pieces. The first one contains
those $\ga$ such that $N\le|\ga|\le\alpha N$ with $\alpha=\,
(\be+\rho)/ (\be-\be')$. This part is bounded by $\alpha N \lambda 
\exp(d\vep\alpha N) $.  It is not difficult to see that if $ \vep d <
\be-\be'$, the second part is also bounded by $\alpha N \lambda
\exp(d\vep\alpha N)$.

\end{pf}

\begin{lem} The following bound holds for $b_2$
  \label{b2}
  \begin{equation}
\label{16}
    b_2 \;\le\; \lambda\, \Bigl({D+1\over \vep}\Bigr)^d
    \left(2\,N\,\Bigl({\be+\rho\over 
    \be-\be^*}\Bigr)\,+\,
 d^d\right)\, \e^{-(\be-\be^*-\vep)N}
  \end{equation}

\end{lem}

\begin{pf}
By Lemma \ref{pab}, $b_2$ is bounded above by 
\begin{eqnarray}
\sum_{\ga\in \G}
(D+|\ga|)^d\, p_\ga\; \sum_{\th\in\G_0} \,\e^{-\be|\th|}
\;+\; \sum_{\ga\in \G}
(D+|\ga|)^d\, \e^{-\be|\ga|}\; \sum_{\th\in\G_0} p_\th
\label{bb9} 
\end{eqnarray}
The first term was controlled in the previous lemma. The second one
can be treated in a similar way.

\end{pf}
\begin{lem}
  \label{b3} The following bound holds for
  $b_3$. For any $\be'\in(\be^*,\be)$, letting $\rho'=\alpha_0(\be')$,
  \begin{equation}
    \label{1b3}
    b_3\;\le\; 2\lambda \,
\e^{-(\be-\be')N}\,+\,Q(\be,\be')\,N\,\lambda D^d\, 
\e^{A(\be,\be') N} \,\e^{-(\be-\be')D} 
  \end{equation}
where $Q(\be,\be')\,=\,
4\,(1-\rho')^{-2}\,(\be-\be')^{-2}\,(2d-1)^2\,(\be+\rho)$  and
$A(\be,\be')\,=\,((\be+\rho)\,(\be-\be')^{-1}+1)(\be + \rho)\log(2d-1) $.
\end{lem}
\begin{pf} Since $\sigma(\eta_0(\th,\C)\,:\, \th \in \G\setminus \B_\ga) 
  \subset \sigma(\A^{\th,0}\,:\, \th \in \G\setminus \B_\ga))$
\begin{eqnarray}
\nonumber
b_3 
 &\le& \sum_{\ga\in \G} \E\Bigl| \E\Bigl\{\eta_0(\ga,\A^\ga) - p_\ga\,|\, 
\sigma(\A^\th\,:\, \th \in \G\setminus \B_\ga)\Bigr\}\Bigr|\\
&=& \sum_{\ga\in \G} \E\Bigl| \E\Bigl\{\eta_0(\ga,\A^\ga) -
\eta_0(\ga,\widehat\A^\ga)|\, 
\sigma(\A^\th\,:\, \th \in \G\setminus \B_\ga)\Bigr\}\Bigr|\\
&\le& \sum_{\ga\in \G} \E \Bigl(\one\Bigl\{\widehat\A^\ga \not\sim
\bigcup_{\th \in \G\setminus \B_\ga} \widehat\A^\th\Bigr\}\,
(\eta_0(\ga,\A^\ga)+\eta_0(\ga,\widehat\A^\ga))\Bigr) \label{bb3}
\end{eqnarray}
where the identity comes from property (c) of the coupling defined in
Section 2 and the second inequality from property (d).
We divide the sum in \reff{bb3} in two parts. The first one is over
those $\ga\in \G$ of length bigger than $KN$, for some $K$ to be
chosen later. Dominating the indicator function by one, this sum is
bounded by
\begin{eqnarray}
  \sum_{\ga\in \G,\,|\ga|\ge KN} 2 p_\ga
&\le& {2 \lambda}\,\e^{(\be+\rho)N}\,\e^{-(\be - \be^*)KN} \label{p18}
\end{eqnarray}
The second sum is over those $\ga\in \G$ with
length less than $KN$. Dominating $\eta_0(\cdot)$ by one, the summand is
bounded by four times
\begin{eqnarray}
 \P\Bigl(\widehat\A^\ga \not\sim
\bigcup_{\th \in \G\setminus \B_\ga(N)} \widehat\A^\th\Bigr) &\le&
(1-\rho')^{-2}\, \sum_{x\in\ga}\sum_{y:d(\ga,\{y\})>D}
|x-y|\, \e^{-(\be-\be')|x-y|} \nonumber\\
&\le& {|\ga|\,D^d\over (1-\rho')^2(\be-\be')}\, \e^{-(\be-\be')D}
\label{p20} 
\end{eqnarray}
using \reff{com} in the first inequality.
Hence the second sum is bounded by
\begin{eqnarray}
&& {{4 \lambda}\,(KN)\,(2d-1)^{NK+2}\,D^d 
  \over(1-\rho')^{2}\,(\be-\be')} 
 \e^{-(\be-\be')D}\,\e^{(\be+\rho)N}\,
\label{p21} 
\end{eqnarray}
using \reff{100} to bound $|V|$. The sum \reff{bb3} is then bounded
above by the sum of \reff{p18} and \reff{p21}. Fixing $K=
{\be+\rho\over\be-\be'}+1$, we get \reff{1b3}.  
\end{pf}
\begin{pf*}{Proof of Theorem \ref{191}.} It remains to choose $D$ in
such a way that $2(2b_1+2b_2+b_3)$ is bounded above by the right hand side
of \reff{125}. We already know that $b_1$ and $b_2$ decay as
$\e^{-(\be-\be^*-\vep)N}$. If we choose $D=\delta N$ in \reff{1b3}, where
\begin{equation}
  \label{dd1}
  \delta\;=\; {\be+\rho\over\be-\be'}\,\Bigl[1+\log(2d-1)\Bigl({1\over
  \be-\be^*}+{1\over\be+\rho} \Bigr)\Bigr] \,+\,1\,,
\end{equation}
then
\begin{equation}
  \label{dd2}
  b_3\;\le\;{4\, \lambda\,(\be+\rho)^d\,N^{d+1}  \over (1-\rho')^2
  (\be-\be')^{2d+2}}\, 
  \, \e^{-(\be-\be')N}\,. 
\end{equation}
With this bound and the choice $\vep=\be'-\be^*$ in \reff{2b2} and
\reff{16} we get \reff{125}.  The constant $M(\be,\be')$ is given by
\begin{equation}
  \label{mbe}
  M(\be,\be')\,=\,10\, 
  {(\be+\rho)^d\over(\be-\be')^{2d+2}}\,\left({d^d\over
  (\be'-\be^*)^{d}}\,+\,{2(\log (2d-1))^2 \over   (1-\rho')^{2}}\, \right)
\end{equation}
\end{pf*}

\begin{ack}
We thank Leonardo Moledo for his contribution to Lemma \ref{pab} and
Roberto Fern\'andez for many discussions, specially about plaquettes.

This paper started when {\rm PAF} was visiting CPT at Marseille. It
was written when PAF was Charg\'e de Recherche associ\'e at
Universit\'e de Rouen. PAF thanks support of Funda\c c\~ao de Amparo
\`a Pesquisa do Estado de S\~ao Paulo and FINEP. PP used financial
support of the project Cofecub-USP 45/97.
\end{ack}

\end{document}